\documentclass[a4paper,12pt]{article}
\usepackage{a4wide}
\usepackage{amsmath}
\usepackage{amssymb}
\usepackage{amsthm}
\usepackage{latexsym}
\usepackage{graphicx}
\usepackage[english]{babel}
\usepackage{makeidx}

\newtheorem{dfn} [subsection]{Definition}
\newtheorem{obs} [subsection]{Remark}
\newtheorem{exm} [subsection]{Example}

\newtheorem{prop}[subsection]{Proposition}

\newtheorem{teor}[subsection]{Theorem}
\newtheorem{lema}[subsection]{Lemma}
\newtheorem{cor} [subsection]{Corollary}

\newcommand{\Nn}{\mathbb N^{n}}

\newcommand{\de}{\mathbf{d}}
\newcommand{\me}{\mathbf{m}}
\newcommand{\leqd}{\leq_{\mathbf{d}}}

\begin{document}

\selectlanguage{english}
\frenchspacing

\large
\begin{center}
\textbf{A generalization of Pardue's formula.}

\normalsize
by
\large

Mircea Cimpoea\c s
\end{center}

\normalsize

\footnotetext[1]{This paper was supported by the CEEX Program of the Romanian
Ministry of Education and Research, Contract CEX05-D11-11/2005.}

\begin{abstract}
In this paper, we introduce a new class of monomial ideals, called $\de$-fixed ideals, which generalize the class of $p$-Borel ideals and show how some results for $p$-Borel ideals can be transfered to this new class. In particular, we give the form of a principal $\de$-fixed ideal and we compute the socle of factors of this ideals, using
methods similar as in \cite{hpv}. This allowed us to give a generalization of Pardue's
formula, i.e. a formula of the regularity for a principal $\de$-fixed ideal.

\vspace{5 pt} \noindent \textbf{Keywords:} p-Borel ideals, Betti
numbers, Mumford-Castelnuovo regularity.

\vspace{5 pt} \noindent \textbf{2000 Mathematics Subject
Classification:}Primary: 13P10, Secondary: 13D25, 13D02, 13H10.
\end{abstract}

\begin{center}
\textbf{Introduction.}
\end{center}

A $p$-Borel ideal is a monomial ideal which satisfy certain combinatorial condition, where $p>0$ is a prime number.
It is well known that any positive integer $a$ has an unique $p$-adic decomposition $a=\sum_{i\geq 0}a_{i}p^{i}$.
 If $a,b$ are two positive integers, we write $a\leq_{p} b$ iff $a_{i}\leq b_{i}$ for any $i$,
 where $a=\sum_{i\geq 0}a_{i}p^{i}$ and $b=\sum_{i\geq 0}b_{i}p^{i}$. We say that a monomial ideal
 $I\subset S=k[x_{1},\ldots,x_{n}]$ is $p$-Borel if for any monomial $u\in I$ and for any indices $j<i$,
 if $t \leq_{p} \nu_{i}(u)$ then $x_{j}^{t}u/x_{i}^{t} \in I$, where $\nu_{i}(u) = max\{k:\; x_{i}^{k}|u\}$.

This definition suggest a natural generalization. The idea is to consider a strictly increasing sequence of positive integers $\de: 1=d_{0}|d_{1}|\cdots|d_{s}$, which we called a $\de$-sequence. Lemma $1.1$ states that
for any positive integer $a$, there exists an unique decomposition
$a=\sum_{t= 0}^{s}a_{t}d_{t}$. If $a,b$ are two positive integers, we write $a\leq_{\de} b$ iff $a_{i}\leq b_{i}$ for any $i$, where $a=\sum_{i\geq 0}a_{i}d_{i}$ and $b=\sum_{i\geq 0}b_{i}d_{i}$. We say that a monomial ideal $I$ is $\de$-fixed if for any monomial $u\in I$ and for any indices $j<i$, if $t \leq_{\de} \nu_{i}(u)$ then $x_{j}^{t}u/x_{i}^{t} \in I$. Obvious, the $p$-Borel ideals are a special case of $\de$-fixed ideals for $\de: 1|p|p^{2}|\cdots$.

A principal $\de$-fixed ideal, is the smallest $\de$-fixed ideal
which contain a given monomial. $1.6$ and $1.8$ gives the explicit
form of a principal $\de$-fixed ideal. In the second section we
compute the socle of factors for a principal $\de$-fixed ideal
($2.1$ and $2.4$). The proofs are similar as in \cite{hpv} but we
consider that is necessary to present them in this context. In the
third section we give a formula ($3.1$) for the regularity of a
$\de$-fixed ideal, which generalize the Pardue's formula for the
regularity of principal $p$-Borel ideals, proved by Aramova-Herzog
\cite{ah} and Herzog-Popescu \cite{hp}. Using a theorem of Popescu
\cite{p} we compute the extremal Betti numbers of $S/I$ ($3.3$).
Also, we show that if $I$ is a principal $\de$-fixed ideal generated
by the power of a variable, then $I_{\geq e}$ is stable for any
$e\geq reg(I)$ ($3.6$), so $reg(I)=min\{e\geq deg(I):I_{\geq e} \ is
\ stable\}$ (3.9). Thus a result of Eisenbud-Reeves-Totaro
\cite[Proposition 12]{ert} holds also in this frame.

The author wish to thanks to his Ph.D.advisor, Professor Dorin
Popescu, for support, encouragement and valuables observations on
the contents of this paper.

\vspace{9mm}

\section{$\de$-fixed ideals.}

In the following $\de: 1=d_{0}|d_{1}|\cdots|d_{s}$ is a strictly increasing sequence of positive integers. We say that
$\de$ is a $\de$-sequence.

\begin{lema}
Let $\de$ be a $\de$-sequence. Then, for any $a\in \mathbb N$, there exists an unique sequence of positive integers $a_{0},a_{1},\ldots,a_{s}$ such that:
\begin{enumerate}
    \item $a= \sum_{t=0}^{s}a_{t}d_{t}$ and
    \item $0 \leq a_{t} < \frac{d_{t+1}}{d_{t}}$, for any $0\leq t < s$.
\end{enumerate}
Conversely, if $\de: 1=d_{0}<d_{1}<\cdots<d_{s}$ is a sequence of positive integers such that for any
$a \in \mathbb N$ there exists an unique sequence of positive integers $a_{0},a_{1},\ldots,a_{s}$ as before, then
$\de$ is a $\de$-sequence.
\end{lema}

\begin{proof}
Let $a_{s}$ be the quotient of $a$ divided by $d_{s}$. For $0\leq t<s$ let $a_{t}$ be the quotient of $(a-q_{t+1})$ divided by $d_{t}$, where $q_{t+1}=\sum_{j=t+1}^{s}a_{j}d_{j}$.
We will prove that $a_{0},a_{1},\ldots,a_{s}$ fulfill the required conditions. Indeed, it is obvious that $a= \sum_{t=0}^{s}a_{t}d_{t}$. On the other hand, $a-q_{t+1}<d_{t+1}$, since $a-q_{t+1}$ is $a-q_{t+2}$ modulo $d_{t+1}$. Therefore, since $a_{t}$ is the quotient of $(a-q_{t+1})$ divided by $d_{t}$, it follows that $a_{t}<\frac{d_{t+1}}{d_{t}}$.

Suppose there exists another decomposition $a = \sum_{j=0}^{s}b_{j}d_{j}$ which also fulfill the conditions $1$ and $2$. Then, we may assume that there exists an integer $0\leq t\leq s$ such that $b_{s}=a_{s}, \cdots, b_{t+1}=a_{t+1}$ and $b_{t}>a_{t}$.
Notice that $d_{t}>\sum_{j=0}^{t-1}a_{j}d_{j}$. Indeed, \[ \sum_{j=0}^{t-1}a_{j}d_{j}\leq \sum_{j=0}^{t-1}(\frac{d_{j+1}}{d_{j}} - 1)d_{j} = (d_{1}-d_{0}) + (d_{2}-d_{1}) + \cdots + (d_{t}-d_{t-1}) = d_{t} - 1<d_{t}. \]
We have $0 = \sum_{j=0}^{s}(b_{j}-a_{j})d_{j} = \sum_{j=0}^{t}(b_{j}-a_{j})d_{j}$, but on the other hand:
\[ (b_{t}-a_{t})d_{t} \geq d_{t} > \sum_{j=0}^{t-1}a_{j}d_{j} \geq \sum_{j=0}^{t-1}(a_{j}-b_{j})d_{j} \]
and therefore $(b_{t}-a_{t})d_{t} - \sum_{j=0}^{t-1}(a_{j}-b_{j})d_{j} = \sum_{j=0}^{t}(b_{j}-a_{j})d_{j} > 0$, which is a contradiction.

For the converse, we use induction on $0 \leq t< s$, the assertion being obvious for $t=0$. Suppose $t>0$ and $d_{0}|d_{1}|\cdots |d_{t}$ and consider the decomposition of $d_{t+1}-1$. Since $d_{t+1}-1<d_{t+1}$, it follows that $d_{t+1}-1 = \sum_{j=0}^{t}a_{t}d_{t}$. On the other hand, since $d_{t+1}-1$ is the largest integer less than  $d_{t+1}$, each $a_{j}$ is maximal between the integers $<d_{j+1}/d_{j}$, for $j<t$. Therefore $a_{j}=d_{j+1}/d_{j}-1$ for $0 \leq j< t$. Thus:
\[ d_{t+1}  = 1 +d_{t+1} - 1  = 1 + a_{0}d_{0} + a_{1}d_{1} + \cdots + a_{t}d_{t} = d_{1} + a_{1}d_{1} + a_{2}d_{2} + \cdots + a_{t}d_{t} = \]\[ = d_{2} + a_{2}d_{2} + \cdots + a_{t}d_{t} = \cdots = (a_{t}+1)d_{t},\;\;so\;\; d_{t}|d_{t+1}.\]
\end{proof}

\begin{dfn}
Let $a,b$ be two positive integers and consider the $\de$-decompositions \linebreak $a=\sum_{j=0}^{s}a_{j}d_{j}$ and $b=\sum_{j=0}^{s}b_{j}d_{j}$. We say that $a\leq_{\de}b$ if $a_{i}\leq b_{i}$ for any $0\leq i\leq s$.
\end{dfn}

\begin{lema}
Let $a,b$ be two positive integers with $a\leq_{\de} b$. Suppose $b=b'+b''$, where $b'$ and $b''$ are positive integers. Then, there exists some positive integers $a'\leqd b'$ and $a''\leqd b''$ such that $a=a'+a''$.
\end{lema}

\begin{proof}
Let $a=\sum_{t=0}^{s}a_{t}d_{t},\; b = \sum_{t=0}^{s}b_{t}d_{t}, \; b'=\sum_{t=0}^{s}b'_{t}d_{t}, \; b''=\sum_{t=0}^{s}b''_{t}d_{t}$. The hypothesis implies $a_{t}\leq b_{t} < d_{t+1}/d_{t}$ and $b'_{t}, b''_{t} < d_{t+1}/d_{t}$ for any $0\leq t<s$. We construct the sequences $a'_{t}, a''_{t}$ using decreasing induction on $t$.
Suppose we have already defined $a'_{j}, a''_{j}$ for $j>t$ such that $\sum_{i=j}^{s}(a'_{i}+a''_{i})d_{i} = \sum_{i=j}^{s}a_{i}d_{i}$ and $b_{t+1} = b'_{t+1}+b''_{t+1}$. This is obvious for $t=s$.

We consider two cases. If $b_{t} = b'_{t} + b''_{t}$, then we choose $a'_{t}\leq b'_{t}$ and $a''_{t}\leq b''_{t}$ such that $a'_{t}+a''_{t}=a_{t}$. We can do this, because $a_{t}\leq b_{t}$. Also, it is obvious from the induction hypothesis that $\sum_{i=t}^{s}(a'_{i}+a''_{i})d_{i} = \sum_{i=t}^{s}a_{i}d_{i}$, so we can pass from $t$ to $t-1$.

If $b_{t} \neq b'_{t} + b''_{t}$ we claim that $b'_{t} + b''_{t} = b_{t}-1$. Indeed, $\sum_{j=0}^{t-1}(b'_{j}+b''_{j})d_{j}<2d_{t}$ and therefore it is impossible to have $b'_{t}+b''_{t} \leq b_{t}-2$, otherwise $\sum_{j=0}^{t}(b'_{j}+b''_{j})d_{j}<b_{t}d_{t}$ and we contradict the equality $b=b'+b''$. Also, since $b'_{t+1} + b''_{t+1} = b_{t+1}$, we cannot have $b'_{t}+b''_{t}>b_{t}$. Similarly we get $b'_{t-1} + b''_{t-1} > b_{t-1}$. By recurrence, we conclude that there exists an integer $u<t$ such that:
$ b'_{u-1}+b''_{u-1} = b_{u-1},\;,b'_{u}+b''_{u} = b_{u} + d_{u+1}/d_{u}$,
$b'_{u+1} + b''_{u+1} = b_{u+1} + d_{u+2}/d_{u+1} - 1,\; \ldots, \; b'_{t-1} + b''_{t-1} = b_{t-1} + d_{t}/d_{t-1} - 1$.

If $a_{j}=b_{j}$ for any $j\in \{u,\ldots,t\}$, we simply choose $a'_{j}=b'_{j}$ and $a''_{j}=b''_{j}$ for any
$j\in \{u,\ldots,t\}$ and the required conditions are fulfilled, so we can pass from $t$ to $u-1$. If this is not the case, then there exists an integer $u \leq q \leq t$ such that $a_{t}=b_{t},\ldots,a_{q+1} = b_{q+1}$ and $a_{q}<b_{q}$. If $q=t$ then for any $j\in \{u,\ldots,t\}$ we can choose $a'_{j}\leq b'_{j}$ and $a''_{j}\leq b''_{j}$ such that $a'_{j}+a''_{j} = a_{j}$. For $j<t$ the previous assertion is obvious because $b'_{j}+b''_{j}\geq b_{j}$, and for $j=t$, since $a_{t}<b_{t}$ we have in fact $a_{t}\leq b'_{t}+b''_{t} = b_{t}-1$ and therefore we can choose again $a'_{t}$ and $a''_{t}$. The conditions are satisfied so we can pass from $t$ to $u-1$.

Suppose $q<t$. For $j\in \{u,\ldots,q-1\}$ we choose $a'_{j}\leq b'_{j}$ and $a''_{j}\leq b''_{j}$ such that $a'_{j}+a''_{j}=a_{j}$. We can do this because $b'_{j}+b''_{j}\geq b_{j} \geq a_{j}$. We choose
$a'_{q}$ and $a''_{q}$ such that $a'_{q}+a''_{q} =a_{q} + d_{q+1}/d_{q}$. We can make this choice, because $a_{q}\leq b_{q}-1$ and $b'_{q}+b''_{q} \geq b_{q} + d_{q+1}/d_{q} - 1$. For $j>q$, we simply put $a'_{j}=b'_{j}$ and $a''_{j}=b''_{j}$. To pass from $t$ to $u-1$ is enough to see that $\sum_{j=u}^{t}a_{j}d_{j} = \sum_{j=u}^{t}(a'_{j}+a''_{j})d_{j}$. Indeed,
\[ \sum_{j=u}^{t}(a'_{j}+a''_{j})d_{j} = \sum_{j=u}^{q-1}(a'_{j}+a''_{j})d_{j} + (a'_{q}+a''_{q})d_{q} + \sum_{j=q+1}^{t}(a'_{j}+a''_{j})d_{j} = \]
\[ = \sum_{j=u}^{q-1}a_{j}d_{j} + (a_{q}+d_{q+1}/d_{q})d_{q} + \sum_{j=q+1}^{t-1}(a_{j}+ d_{j+1}/d_{j}-1)d_{j} + (a_{t}-1)d_{t} =\]\[ = \sum_{j=u}^{t}a_{j}d_{j} + d_{q+1} + \sum_{j=q+1}^{t-1}(d_{j+1}-d_{j}) - d_{t} = \sum_{j=u}^{t}a_{j}d_{j},\]
The induction ends when $t=-1$. Finally, we obtain $a'$ and $a''$ such that $a'+a''=a$, $a'_{t}\leq b'_{t}$ and $a''_{t}\leq b''_{t}$, as required.
\end{proof}

\begin{dfn}
We say that a monomial ideal $I\subset S = k[x_{1},\ldots,x_{n}]$ is $\de$-fixed, if for any monomial $u\in I$ and for any indices $1\leq j<i \leq n$, if $t\leq_{\de} \nu_{i}(u)$ (where $\nu_{i}(u)$ denotes the exponent of the variable $x_{i}$ in $u$) then $u \cdot x_{j}^{t}/x_{i}^{t} \in I$.

Notice that if $\de: 1|p|p^{2}|p^{3}|\cdots$ then $I$ is a $p$-Borel ideal.
\end{dfn}

\begin{dfn}
A $\de$-fixed ideal $I$ is called principal if it is generated, as a $\de$-fixed ideal by one monomial $u$, i.e. $I$ is the smallest $\de$-fixed ideal which contain $u$. We write $I = <u>_{\de}$.
More generally, if $u_{1},\ldots,u_{r} \in S$ are monomials, the $\de$-fixed ideal generated by $u_{1},\ldots,u_{r}$ is the smallest $\de$-fixed ideal $I$ which contains $u_{1},\ldots,u_{r}$. We write $I = <u_{1},\ldots,u_{r}>_{\de}$.
\end{dfn}

Our next goal is to describe the principal $\de$-fixed ideals. The easiest case is when we have a $\de$-fixed ideal generated by the power of a variable. Denote $\me=(x_{1},\ldots,x_{n})$ and $\me^{[d]}=(x_{1}^{d},\ldots,x_{n}^{d})$ for some nonnegative integer $d$. We have the following proposition.

\begin{prop}
If $u=x_{n}^{\alpha}$, then $I= <u>_{\de} = \prod_{t=0}^{s} (\me^{[d_{t}]})^{\alpha_{t}}$,
where $\alpha=\sum_{t=0}^{s}\alpha_{t}d_{t}$.
\end{prop}

\begin{proof}
Let $I' = \prod_{t=0}^{s} (\me^{[d_{t}]})^{\alpha_{t}}$. The minimal generators of $I'$ are monomials of the type
$w = \prod_{t=0}^{s}\prod_{j=1}^{n} x_{j}^{\lambda_{tj}\cdot d_{t}}$, where $0\leq \lambda_{tj}$ and $\sum_{j=1}^{n}\lambda_{tj} = \alpha_{t}$. First, let us show that $I'\subset I$. In order to do this, we choose $w$
a minimal generator of $I'$ (the one bellow). We write $x_{n}^{\alpha}$ like this: $x_{n}^{\alpha} = x_{n}^{\alpha_{0}d_{0} + \alpha_{1}d_{1} + \cdots + \alpha_{s}d_{s}} = x_{n}^{\alpha_{0}d_{0}} \cdot x_{n}^{\alpha_{1}d_{1}} \cdots x_{n}^{\alpha_{s}d_{s}}$.
Since $\lambda_{01}d_{0}\leq_{\de}\alpha_{0}d_{0} + \alpha_{1}d_{1} + \cdots + \alpha_{s}d_{s}$ and $I$ is $\de$-fixed it follows that $x_{1}^{\lambda_{01}d_{0}}x_{n}^{\alpha-\lambda_{01}d_{0}} \in I$. Also, $\lambda_{02}d_{0} < \alpha-\lambda_{01}d_{0} = (\alpha_{0}-\lambda_{01})d_{0} + \alpha_{1}d_{1} + \cdots + \alpha_{s}d_{s}$, and since $I$ is $\de$-fixed it follows that $x_{1}^{\lambda_{01}d_{0}} x_{2}^{\lambda_{02}d_{0}}x_{n}^{\alpha-\lambda_{01}d_{0}-\lambda_{02}d_{0}} \in I$.
Using iteratively this argument, one can easily see that $x_{1}^{\lambda_{01}d_{0}} \cdots x_{n}^{\lambda_{0n}d_{0}}x_{n}^{\alpha-\alpha_{0}d_{0}} \in I$. Also $\alpha - \alpha_{0}d_{0} = \alpha_{1}d_{1} + \cdots +\alpha_{n}d_{n}$. Again, using an inductive argument, we get:
\[ (x_{1}^{\lambda_{01}d_{0}} \cdots x_{n}^{\lambda_{0n}d_{0}})\cdot (x_{1}^{\lambda_{11}d_{1}} \cdots x_{n}^{\lambda_{1n}d_{1}})\cdots (x_{1}^{\lambda_{s1}d_{s}} \cdots x_{n}^{\lambda_{sn}d_{s}}) = w \in I.\]
For the converse, i.e. $I\subset I'$, is enough to verify that $I'$ is $\de$-fixed. In order to do this, is enough to
prove that the minimal generators of $I'$ fulfill the definition of a $\de$-fixed ideal. Let $w$ be a minimal generator of  $I'$. Let $2\leq i \leq n$. Then $\nu_{i}(w) = \sum_{t=0}^{s}\lambda_{ti}d_{t}$. If $\beta\leqd \nu_{i}(w)$ then  $\beta=\sum_{t=0}^{s}\beta_{t}d_{t}$ with $\beta_{t}\leq \lambda_{ti}$. Let $1\leq k<i$. We have
 \[ w\cdot x_{k}^{\beta}/x_{i}^{\beta} = \prod_{t=0}^{s}(\prod_{j\neq i,k}x_{j}^{\lambda_{tj}d_{t}})\cdot x_{i}^{(\lambda_{ti} - \beta_{t})d_{t}}\cdot x_{k}^{(\lambda_{tk} + \beta_{t})d_{t}}. \]
Thus $w\cdot x_{k}^{\beta}/x_{i}^{\beta} \in I'$ and therefore $I'$ is $\de$-fixed. Since $I$ is the smallest $\de$-fixed ideal which contains $x_{n}^{\alpha}$ it follows that $I\subset I'$.
\end{proof}

\begin{prop}
If $\alpha \leq \beta$ then $<x_{n}^{\beta}>_{\de} \subseteq <x_{n}^{\alpha}>_{\de}$.
\end{prop}

\begin{proof}
The case $\alpha=\beta$ is obvious, so we may assume $\alpha<\beta$. We denote $I=<x_{n}^{\alpha}>_{\de}$ and $I' = <x_{n}^{\beta}>_{\de}$. We write $\alpha=\sum_{t=0}^{s}\alpha_{t}d_{t}$ and $\beta=\sum_{t=0}^{s}\beta_{t}d_{t}$. If $w$ is a minimal generator of $I'$ then $w=\prod_{t=0}^{s}\prod_{i=1}^{n}x_{i}^{\lambda_{ti}d_{t}}$, where $0\leq \lambda_{ti}$ and $\sum_{i=1}^{n}\lambda_{ti}=\beta_{t}$. We claim that $w\in I$ and therefore $I'\subset I$ as required.

Since $\alpha < \beta$ there exists $t\in \{0,\ldots,s\}$ such that $\alpha_{s} = \beta_{s},\ldots,\alpha_{t+1}=\beta_{t+1}$ and $\alpha_{t}<\beta_{t}$. We may assume $\lambda_{t1}>0$. We have
\[ w = \prod_{j=0}^{s}\prod_{i=1}^{n}x_{i}^{\lambda_{ji}d_{j}} = \prod_{j=0}^{t-1}x_{1}^{\alpha_{j}d_{j}} x_{1}^{(\lambda_{t1}-1)d_{t}} x_{1}^{d_{t}-\sum_{j=0}^{t-1}\alpha_{j}d_{j}} \prod_{i=2}^{n}x_{i}^{\lambda_{ti}d_{t}}\prod_{j>t}\prod_{i=1}^{n}x_{i}^{\lambda_{ji}d_{j}} \]
and now it is obvious that $w\in I$.
\end{proof}

We have the general description of a principal $\de$-fixed ideal given by the following proposition. In the proof, we will apply Lemma $1.3$.

\begin{prop}
Let $1\leq i_{1}<i_{2}<\cdots<i_{r}=n$ and let $\alpha_{1},\ldots,\alpha_{r}$ be some positive integers.
If $u=x_{i_{1}}^{\alpha_{1}}x_{i_{2}}^{\alpha_{2}} \cdots x_{i_{r}}^{\alpha_{r}}$ then:
\[ I=<u>_{\de} = <x_{i_{1}}^{\alpha_{1}}>_{\de}\cdot <x_{i_{2}}^{\alpha_{2}}>_{\de} \cdots <x_{i_{r}}^{\alpha_{r}}>_{\de} = \prod_{q=1}^{r}\prod_{t=0}^{s}(\me_{q}^{[d_{t}]})^{\alpha_{qt}},\]
where $\me_{q} = (x_{1},x_{2},\ldots,x_{i_{q}})$ and $\alpha_{q}=\sum_{t=0}^{s}\alpha_{qt}d_{t}$.
\end{prop}

\begin{proof}
Let $I' = \prod_{q=1}^{r}\prod_{t=0}^{s}(\me_{q}^{[d_{t}]})^{\alpha_{qt}}$. The minimal generators of $I'$ are monomials of the type $w = \prod_{q=1}^{r} \prod_{t=0}^{s}\prod_{j=1}^{i_{q}} x_{j}^{\lambda_{qtj}\cdot d_{t}}$, where $0\leq \lambda_{qtj}$ and $\sum_{j=1}^{n}\lambda_{qtj} = \alpha_{qt}$. First, we show that $I'\subset I$. In order to do this,
it is enough to prove that by iterative transformations we can modify $u$ such that we obtain $w$.

The idea of this transformations is the same as in the proof of $1.6$. Without given all the details, one can see that if we rewrite $u$ as \[ (x_{i_{1}}^{\alpha_{10}d_{0}}x_{i_{1}}^{\alpha_{11}d_{1}}\cdots x_{i_{1}}^{\alpha_{1s}d_{s}}) \cdots (x_{i_{r}}^{\alpha_{r0}d_{0}}x_{i_{r}}^{\alpha_{r1}d_{1}}\cdots x_{i_{r}}^{\alpha_{rs}d_{s}} ),\]
where $\alpha_{q}=\sum_{t=0}^{s}\alpha_{qt}d_{t}$, we can pass to $w$, using the transformations
\[x_{i_{1}}^{\alpha_{10}d_{0}} \mapsto \prod_{j=1}^{i_{1}}x_{j}^{\lambda_{10j}d_{0}},\ldots, x_{i_{1}}^{\alpha_{1s}d_{s}} \mapsto \prod_{j=1}^{i_{1}}x_{j}^{\lambda_{1sj}d_{s}}, \ldots , x_{i_{r}}^{\alpha_{r0}d_{0}} \mapsto \prod_{j=1}^{i_{r}}x_{j}^{\lambda_{r0j}d_{0}}, \ldots, x_{i_{r}}^{\alpha_{rs}d_{s}} \mapsto \prod_{j=1}^{i_{r}}x_{j}^{\lambda_{rsj}d_{s}}.\]
Therefore $w \in I$, and thus $I'\subset I$. For the converse, it is enough to see that $I'$ is a $\de$-fixed ideal. Let $w$ be a minimal generator of $I'$. We choose an index $2 \leq i \leq n$. Then $\nu_{i}(w) = \sum_{q=1}^{r}\sum_{t=0}^{s}\lambda_{qti}d_{t}$. Let $\beta \leq \nu_{i}(w)$. Using Lemma $3.1$, we can choose some positive integers $\beta_{1},\ldots,\beta_{r}$ such that:
\[ (a) \beta = \sum_{q=1 , i_{q}\geq i}^{r} \beta_{q} \;and\; (b) \beta_{q}\leqd \sum_{t=0}^{s}\lambda_{qti}d_{t},\; \]
i.e. $\beta_{qt}\leq \lambda_{qti}$, where $\beta_{q} = \sum_{t=0}^{s}\beta_{qt}d_{t}$. Let $k<i$. Then,
\[ w\cdot x_{k}^{\beta}/x_{i}^{\beta} = \prod_{q=1}^{r}  \prod_{t=0}^{s}\left( \prod_{j=1, j\neq k,i}^{i_{q}} x_{j}^{\lambda_{qtj}\cdot d_{t}} \right)
x_{i}^{(\lambda_{qti} - \beta_{qt})d_{t}} x_{k}^{(\lambda_{qtk} + \beta_{qt})d_{t}}.\]
Now, it is easy to see that $w\cdot x_{k}^{\beta}/x_{i}^{\beta} \in I'$, and therefore $I'$
is $\de$-fixed.
\end{proof}

\begin{exm}
Let $\de: 1|2|4|12$.
\begin{enumerate}
    \item Let $u=x_{3}^{21}$. We have $21=1\cdot 1 + 0\cdot 2 + 2 \cdot 4 + 1 \cdot 12$. From $1.6$, we get:
    \[ <u>_{\de} = (x_{1},x_{2},x_{3})(x_{1}^{4},x_{2}^{4},x_{3}^{4})^{2}(x_{1}^{12},x_{2}^{12},x_{3}^{12}).\]
    \item Let $u=x_{1}^{2}x_{2}^{9}x_{3}^{16}$. We have $9=1\cdot 1 + 2\cdot 4$ and $16 = 1\cdot 4 + 1\cdot 12$. From
          $1.8$, we get
    \[ <u>_{\de} = x_{1}^{2}<x_{2}^{9}>_{\de}<x_{3}^{16}>_{\de} = x_{1}^{2}(x_{1},x_{2})(x_{1}^{4},x_{2}^{4})^{2}
                    (x_{1}^{4},x_{2}^{4},x_{3}^{4})(x_{1}^{12},x_{2}^{12},x_{3}^{12}).\]
\end{enumerate}
\end{exm}

\begin{dfn}
We say that a monomial ideal $I\subset k[x_{1},\ldots,x_{n}]$ is a Borel type ideal if
\[ I:x_{j}^{\infty} = I:(x_{1},\ldots,x_{j})^{\infty},\;for\;any\;j=1,\ldots,n.\]
\end{dfn}

\begin{prop}
Any $\de$-fixed ideal $I$ is a Borel type ideal.
\end{prop}

\begin{proof}
Indeed, $[3,\; Proposition \; 2.2]$ says that an ideal $I$ is of
Borel type if and only if for any $1\leq j<i\leq n$, there exists an
positive integer $t$ such that $x_{j}^{t}(u/x_{i}^{\nu_{i}(u)})\in
I$. Choosing $t=\nu_{i}(u)$, is easy to see that the definition of a
$\de$-fixed ideal implies the condition above.
\end{proof}

\begin{dfn}
Let $S=k[x_{1},\ldots,x_{n}]$ and let $M$ be a finitely generated graded $S$-module. The module $M$ is sequentially Cohen-Macaulay if there exists a finite filtration $0 = M_{0} \subset M_{1} \subset \cdots \subset M_{r} = M$ of $M$
by graded submodules of $M$ such that:
\begin{itemize}
    \item $M_{i}/M_{i-1}$ are Cohen-Macaulay for any $i=1,\ldots,r$ and
    \item $dim(M_{1}/M_{0})< dim(M_{2}/M_{1})< \cdots < dim(M_{r}/M_{r-1})$.
\end{itemize}
In particular, if $I\subset S$ is a graded ideal then $R=S/I$ is sequentially Cohen-Macaulay if there exists a chain of ideals $I = I_{0} \subset I_{1} \subset \cdots \subset I_{r} = S$ such that $I_{j}/I_{j-1}$ are Cohen-Macaulay and $dim(I_{j}/I_{j-1})<dim(I_{j+1}/I_{j})$ for any $j=1,\ldots,r-1$.
\end{dfn}

\begin{obs}
Let $I\subset S$ be a monomial ideal. Recursively we define an ascending chain of monomial ideals as follows: We let $I_{0}:=I$. Suppose $I_{\ell}$ is already defined. If $I_{\ell}=S$ then the chain ends. Otherwise, let $n_{\ell} = max\{i:\;x_{i}|u$ for an $u\in G(I_{\ell}) \}$. We set $I_{\ell +1}:= (I_{\ell}:x_{n_{\ell}}^{\infty})$. It is obvious that $n_{\ell}>n_{\ell-1}$, and therefore the chain $I_{0}\subset I_{1} \subset \cdots \subset I_{r}=S$ is finite and has length $r\leq n$. We call this chain of ideals, the sequential chain of $I$.

If $I$ is a Borel type ideal, $[3,\;Lemma\;2.4]$ says that $I_{\ell
+1}:= I_{\ell}:(x_{1},x_{2},\ldots,x_{n_{\ell}})^{\infty}$. From
$[3,\;Corollary\;2.5]$, it follows that $R=S/I$ is sequentially
Cohen-Macaulay with the sequential chain $I_{0}\subset I_{1} \subset
\cdots \subset I_{r}=S$ defined above. Moreover $I_{\ell+1}/I_{\ell}
\cong J_{\ell}^{sat}/J_{\ell}[x_{n_{\ell}+1}, \ldots, x_{n}]$, where
$J_{\ell} = I_{\ell}\cap k[x_{1},\ldots,x_{n_{\ell}}]$ and
$J_{\ell}^{sat} = J_{\ell}:(x_{1},\ldots,x_{n_{\ell}})^{\infty}$.

Let $u = x_{i_{1}}^{\alpha_{1}}x_{i_{2}}^{\alpha_{2}} \cdots x_{i_{r}}^{\alpha_{r}}$ and
$ I=<u>_{\de} =  \prod_{q=1}^{r}\prod_{t=0}^{s}(\me_{q}^{[d_{t}]})^{\alpha_{qt}}$,
where $\me_{q} = (x_{1},\ldots,x_{i_{q}})$ and $\alpha_{q}=\sum_{t=0}^{s}\alpha_{qt}d_{t}$.
Let $I_{r-e} = \prod_{q=1}^{e}\prod_{t=0}^{s}(\me_{q}^{[d_{t}]})^{\alpha_{qt}}$. Then $I = I_{0} \subset I_{1} \subset \cdots \subset I_{r}=S$ is the sequential chain of $I$. Let $n_{\ell}=i_{q_{r-\ell}}$.
Indeed, since $x_{n_{\ell}}^{\alpha_{r-e}}I_{\ell+1} \subset I_{\ell} \Rightarrow I_{\ell +1}\subset (I_{\ell}:x_{n_{\ell}}^{\infty})$. For the converse, let $w\in (I_{\ell}:x_{n_{\ell}}^{\infty})$ be any minimal generator. Then there exists an integer $b$ such that $w\cdot x_{n_{\ell}}^{b}\in I_{\ell}$. We may assume that $w$ is a minimal generator of $I_{\ell}$. Then $w\cdot x_{n_{\ell}}^{b} = w'\cdot y$ for a $w'\in I_{\ell+1}$ and $y\in \prod_{j=0}^{t}(\me_{r-\ell}^{[d_{j}]})^{\alpha_{r-\ell,j}}$ with $x_{n_{\ell}}^{b}|y$. Thus $w'|w$, and therefore $w\in I_{\ell+1}$.
\end{obs}

Let $S=k[x_{1},\ldots,x_{n}]$ and let $M$ be a finitely graded generated $S$-module with the minimal graded free resolution $0 \rightarrow F_{s} \rightarrow F_{s-1} \rightarrow \cdots \rightarrow F_{0} \rightarrow M \rightarrow 0$. Let $Syz_{t}(M)=Ker(F_{t} \rightarrow F_{t-1})$. The module $M$ is called $(r,t)$-regular if $Syz_{t}(M)$ is $(r+t)$-regular in the sense that all generators of $F_{j}$ for $t\leq j\leq s$ have degrees $\leq j+r$. The $t$-regularity of $M$ is by definition $(t-reg)(M) = min\{r|\; M\;is\;(r,t)-regular \}$.

Obvious $(t-reg)(M)\leq ((t-1)-reg)(M)$. If the equality is strict and $r=(t-reg)(M)$ then $(r,t)$ is called a corner of $M$ and $\beta_{t,r+t}(M)$ is an extremal Betti number of $M$, where $\beta_{ij} = dim_{k}Tor_{i}(k,M)_{j}$ denotes
the $ij$-th graded Betti number of $M$. Later, we will use the following result:

\begin{teor}\cite[Theorem 3.2]{p}
If $I\subset S$ is a Borel type ideal, then $S/I$ has at most $r+1$-corners among $(n_{\ell},s(J_{\ell}^{sat}/J_{\ell}))$
and the corresponding extremal Betti numbers are
\[ \beta_{n_{\ell},s(J_{\ell}^{sat}/J_{\ell}) + n_{\ell}} (S/I) = dim_{k}(J_{\ell}^{sat}/J_{\ell})_{s(J_{\ell}^{sat}/J_{\ell})}.\]
\end{teor}

\section{Socle of factors by principal $\de$-fixed ideals.}

In the following, we suppose $n\geq 2$.

\begin{lema}
Let $\de:1=d_{0}|d_{1}| \cdots |d_{s}$, $\alpha\in \mathbb N$ and $I=<x_{n}^{\alpha}>_{\de} =\prod_{t=0}^{s} (\me^{[d_{t}]})^{\alpha_{t}}$. Let $q_{t}=\sum_{j=t}^{s}\alpha_{j}d_{j}$. Let
\[ J = \sum_{t=0 , \alpha_{t}>0}^{s} (x_{1}\cdots x_{n})^{d_{t}-1} (\me^{[d_{t}]})^{\alpha_{t}-1} \prod_{j>t} (\me^{[d_{j}]})^{\alpha_{j}}.\]
Then:
\begin{enumerate}
    \item $Soc(S/I) = \frac{J+I}{I}$
    \item Let $e$ be a positive integer. Then $(\frac{J+I}{I})_{e}\neq 0 \Leftrightarrow
          e = q_{t} + (n-1)(d_{t}-1) -1$, for some $0\leq t\leq s$ with $\alpha_{t}>0$.
    \item $max\{e| (\frac{J+I}{I})_{e}\neq 0 \} = \alpha_{s}d_{s} + (n-1)(d_{s}-1) - 1$.
\end{enumerate}
\end{lema}

\begin{proof}
1. First we prove that $\frac{J+I}{I} \subset Soc(S/I)$. Since $Soc(S/I) = (O:_{S/I}\mathbf{m})$, it is enough to show that $\textbf{m}J \subset I$.

We have $J=\sum_{t=0,\;\alpha_{t}>0}^{s}J_{t}$, where $J_{t} = (x_{1}\cdots x_{n})^{d_{t}-1} (\me^{[d_{t}]})^{\alpha_{t}-1} \prod_{j>t} (\me^{[d_{j}]})^{\alpha_{j}}$. It is enough
to prove that $x_{i}J_{t}\subset I$ for any $i$ and any $t$. Suppose $i=1$:
\[x_{1}J_{t} = x_{1}^{d_{t}}(x_{2}\cdots x_{n})^{d_{t}-1} (\me^{[d_{t}]})^{\alpha_{t}-1} \prod_{j>t} (\me^{[d_{j}]})^{\alpha_{j}} \subset (x_{2}\cdots x_{n})^{d_{t}-1} \prod_{j\geq t} (\me^{[d_{j}]})^{\alpha_{j}}. \]
On the other hand, $(x_{2}\cdots x_{n})^{d_{t}-1} \in \prod_{j<t}(\me^{[d_{j}]})^{\alpha_{j}}$, because $d_{t}-1\geq \sum_{j<t}\alpha_{j}d_{j}$. Thus $x_{1}J_{t}\subset I$.

For the converse, we apply induction on $\alpha$. If $\alpha=1$ then $s=0$ and $I=(x_{1},\ldots,x_{n}) = \mathbf{m}$. $J=(x_{1},\ldots,x_{n})^{d_{0}-1}=S$, and obvious $Soc(S/I) = Soc(S/\mathbf{m}) = S/\mathbf{m}$. Let us suppose that
$\alpha>1$. We prove that if $w\in S\setminus I$ is a monomial such that $\mathbf{m}w\subset I$, then $w\in J$. Let $t_{e}=max\{ t:\; x_{e}^{d_{t}-1}|w \}$. Renumbering $x_{1},\ldots,x_{n}$ which does not affect either $I$ or $J$, we may suppose that $t_{1}\geq t_{2}\geq \cdots \geq t_{n}$. We have two cases: (i)$t_{1}>t_{n}$ and (ii)$t_{1}=t_{n}$.
But first, let's make the following remark: $(*)$ If $u=x_{1}^{\beta_{1}}\cdots x_{n}^{\beta_{n}} \in \prod_{j\geq t}\me^{d_{j}}$ and $\beta_{i}<d_{t}$ for certain $i$ then $u/x_{i}^{\beta_{i}} \in \prod_{j\geq t}\me^{[d_{j}]}$ (the proof is similarly to \cite[Lemma 3.5]{hpv}).

In the case (i), there exists an index $e$ such that $t_{e}>t_{e+1} = \cdots = t_{n}$. Then we have $w=(x_{n}\cdots x_{e+1})^{d_{t_{n}}-1}\cdot x_{e}^{d_{t_{e}}-1}\cdot y$, for a monomial $y\in S$. We consider two cases (a) $x_{e}$ does not divide $y$ and (b) $x_{e}$ divide $y$. (a) From $x_{n}w = x_{n}^{d_{t_{n}}}\cdot (x_{n-1}\cdots x_{e+1})^{d_{t_{n}}-1}x_{e}^{d_{t_{e}}-1}\cdot y \in I$ we see that $y \in \prod_{j\geq t_{e}} (\mathbf{m}^{[d_{j}]})^{\alpha_{j}}$, by $(*)$. Therefore $w\in I$, because
$x_{e}^{d_{t_{e}}-1} \in \prod_{j<t_{e}}(\me^{[d_{j}]})^{\alpha_{j}}$, which is an contradiction.

(b) In this case, $w = (x_{n}\cdots x_{e+1})^{d_{t_{n}}-1} x_{e}^{d_{t_{e}}}y'$, where $y'=y/x_{e}$. We claim that there exist $\lambda\leq t_{e}$ such that $\alpha_{\lambda}\neq 0$. Indeed, if all $\alpha_{\lambda}=0$ for $\lambda\leq t_{e}$, then $I = \prod_{j=t_{e}+1}^{s}(\mathbf{m}^{[d_{j}]})^{\alpha_{j}}$ and $x_{n}w\in I$ implies $y'\in I$ because of the maximality of $t_{n}$ and $(*)$. It follows $w\in I$, which is false. Choose $\lambda\leq t_{e}$ maximal possible with $\alpha_{\lambda}\neq 0$. Set $w' = w/ x_{e}^{d_{\lambda}}$. Note that $\mathbf{m}w\subset I$ implies
\[\mathbf{m}w\subset I' = (\mathbf{m}^{[d_{t_{\lambda}}]})^{\alpha_{\lambda}-1} \prod_{j\neq \lambda} (\mathbf{m}^{[d_{t_{j}}]})^{\alpha_{j}}.\]
It is obvious that $x_{q}w'\in I'$ for $q\neq e$. Also, since $x_{e}^{d_{t_{e}+1}}$ does not divide $x_{e}w$ implies $x_{e}w'\in I'$. Choosing $\alpha'=\alpha-d_{\lambda}$, we get $\alpha'_{j}=\alpha_{j}$ for $j\neq \lambda$ and $\alpha'_{\lambda} = \alpha_{\lambda} - 1$ and therefore we can apply our induction hypothesis for $I'$ (because $\alpha'<\alpha$) and for the ideal $J'$ associated to $I'$, which has the form:
\[ J' = \sum_{q=0,\alpha'_{q}\neq 0} (x_{1}\cdots x_{n})^{d_{q}-1} (\mathbf{m}^{[d_{q}]})^{\alpha'_{q}-1} \prod_{j>q} (\mathbf{m}^{[d_{j}]})^{\alpha'_{j}}, \]
and so $w = x_{e}^{d_{\lambda}}w' \in x_{e}^{d_{\lambda}}J' \subset J$.

It remains to consider the case (ii) in which we have in fact $t_{1}=t_{2}=\cdots = t_{n}$. If $y=w/ (x_{1}\cdots x_{n})^{d_{t_{n}}-1} \in \mathbf{m}$, then there exists $e$ such that $x_{e}|y$, and we apply our induction hypothesis as in the case $(b)$ above. Thus we may suppose $y=1$, i.e. $w=(x_{1}\cdots x_{n})^{d_{t_{n}}-1}$. Since $\mathbf{m}w\subset I$, we see that $\alpha_{j}=0$ for $j>t_{n}$ and $\alpha_{t_{n}}=1$ (otherwise $w\in I$, which is absurd). Thus $w\in J$.

2. Let $v=x_{1}^{q_{t}-1}(x_{2}\cdots x_{n})^{d_{t}-1}$. Then $deg(v)=q_{t} + (n-1)(d_{t}-1) -1$. But $v\in J$ and $v\notin I$, therefore $v\neq 0$ in $Soc(S/I) = \frac{J+I}{I}$.

% On the other hand, if $u\neq 0$ \^in $Soc(S/I) = (0:_{S/I}\mathbf{m})$, then $x_{i}u = 0$. Thus, we have nonzero monomials in $Soc(S/I)$ only in degrees $q_{t} + (n-1)(d_{t}-1) -1$.

3. Let $e_{t}= q_{t} + (n-1)(d_{t}-1) -1$ for $0\leq t\leq s$. Let $t<s$. Then
\[ e_{t+1}-e_{t} = q_{t+1}-q_{t} + (n-1)(d_{t+1}-d_{t}) = -\alpha_{t}d_{t} + (n-1)(d_{t+1}-d_{t}) \geq d_{t+1} - (\alpha_{t}+1)d_{t} \geq 0,\;so \]
\[ max\{e| ((J+I)/I)_{e}\neq 0 \} = e_{s} = \alpha_{s}d_{s} + (n-1)(d_{s}-1) - 1. \]
\end{proof}

\begin{obs}
From the proof of the above lemma, we may easily conclude that for
$n\geq 3$, $e_{t}=e_{t'}$ if and only if $t=t'$, and if $n=2$, then $e_{t}=e_{t'}$ ($t<t'$) if and only if \linebreak $\alpha_{t'-1}=d_{t'}/d_{t'-1},\ldots,\alpha_{t}=d_{t+1}/d_{t}$.
\end{obs}

\begin{cor}
With the notations of previous lemma and remark, let $0\leq t\leq s$ be an integer such that $\alpha_{t}\neq 0$. Let $h_{t}=dim_{K}((I+J_{t})/I)$. Then:
\begin{enumerate}
    \item $G(J_{t})\cap (I+J_{t'}) = 0$ for $0\leq t'\leq s$, $t'\neq t$.
    \item $h_{t}=\binom{n+\alpha_{t}-2}{n-1}\prod_{j>t}\binom{n+\alpha_{j}-1}{n-1}$.
    \item $dim_{K}(Soc(S/I)_{e}) =
\begin{cases}
    h_{q},& if\; n\geq 3\;and\;e=e_{q}\;for\;a\;q\leq s\;with\;\alpha_{q}\neq 0. \\
    \sum_{q}h_{q},& if\; n=2\;and\;q\in\{\epsilon|e=e_{\epsilon}\;for\;\epsilon\leq s\;with\;\alpha_{\epsilon}\neq 0\}.\\
    0,& otherwise.
\end{cases}.$
\end{enumerate}
\end{cor}

\begin{proof}
1. First suppose $t'<t$. A minimal generator $x^{\beta}=x_{1}^{\beta_{1}}\cdots x_{n}^{\beta_{n}}$ of $J_{t}$ has the form \[(x_{1}\cdots x_{n})^{d_{t}-1}\prod_{j\geq t}(x_{1}^{\lambda_{1j}d_{j}}\cdots x_{n}^{\lambda_{nj}d_{j}}),\;where \;\sum_{\nu=1}^{n}\lambda_{\nu j} = \begin{cases} \alpha_{j},& if\; j>t,\\ \alpha_{t}-1,&if\; j=t.  \end{cases}.\] Thus, $\beta_{i} = d_{t}-1 + \sum_{j=t}^{s}\lambda_{ij}d_{j}$. On the other hand, $d_{t}-1 = \sum_{j=0}^{t-1} (d_{j+1}/d_{j}-1)d_{j}$, so $\beta_{i}$ has the writing $\sum_{j=0}^{s}\beta_{ij}d_{j}$, where $\beta_{ij}=d_{j+1}/d_{j}-1$ for $j<t$ and $\beta_{ij}=\lambda_{ij}$ for $j\geq t$.

Assume that $x^{\beta}\in I+J_{t'}$ for a certain $t'<t$. Then there exists $\gamma\in\Nn$ such that $x^{\gamma}\in G(I)$ (or $x^{\gamma}\in G(J_{t'})$) and $x^{\gamma}|x^{\beta}$, that is $\gamma_{i}\leq \beta_{i}$ for all $1\leq i\leq n$. Let $\gamma_{i} = \sum_{j=0}^{s}\gamma_{ij}d_{j}$, the $\de-$ decomposition of $\gamma_{i}$. We notice that $(\beta_{is},\ldots,\beta_{i0})\geq (\gamma_{is},\ldots,\gamma_{i0})$ in the lexicographic order.

Note that all minimal generators $x^{\gamma}$ of $I$ have the same degree $\alpha< e_{t}$ and $\sum_{i=1}^{n}\gamma_{iq} = \alpha_{q}$ for each $0\leq q \leq s$. Also
all minimal generators $x^{\gamma}$ of $J_{t'}$ have the same degree $e_{t'}<e_{t}$ and $\sum_{i=1}^{n}\gamma_{iq}=\alpha_{q}$ for each $t\leq q\leq s$. It follows $deg(x^{\beta})>deg(x^{\gamma})$ and so $\beta_{i}>\gamma_{i}$ for some $i$. Choose a maximal $q<s$ such that $\beta_{iq}>\gamma_{iq}$ for some $i$. Thus $\beta_{ij}=\gamma_{ij}$ for $j>q$. It follows $\beta_{iq}\geq \gamma_{iq}$ since $(\beta_{is},\ldots,\beta_{i0})\geq_{lex} (\gamma_{is},\ldots,\gamma_{i0})$. If $q\leq t$ then we have
\[ \alpha_{q} = \sum_{i=1}^{n}\gamma_{iq} < \sum_{i=1}^{n}\beta_{iq} = \sum_{i=1}^{n}\lambda_{iq}\leq \alpha_{q}, \]
which is not possible. It follows $q<t$ and so $\beta_{it} = \gamma_{it}$ for each $i$. But this is not possible because we get $\alpha_{t} = \sum_{i=1}^{n} \gamma_{it} = \sum_{i=1}^{n} \lambda_{it} = \alpha_{t}-1$. Hence $x^{\beta}\notin I+J_{t'}$.

Suppose now $t'>t$. If $e_{t'}>e_{t}$, then $G(J_{t})\cap G(J_{t'}) = \emptyset$ by degree reason. Assume $e_{t}=e_{t'}$. If follows $n=2$ by the previous remark. If $x_{1}^{\beta_{1}}x_{2}^{\beta_{2}}\in G(J_{t})\cap J_{t'}$ we necessarily get
$x_{1}^{\beta_{1}}x_{2}^{\beta_{2}}\in G(J_{t})\cap G(J_{t'})$ again by degree reason. But this is not possible since it implies that $\alpha_{t'}-1 = \beta_{1t'} + \beta_{2t'} = \alpha_{t'}$.

2. and 3. follows from 1.
\end{proof}

\begin{teor}
Let $u = \prod_{q=1}^{r}x_{i_{1}}^{\alpha_{q}}$, where $2\leq i_{1}<i_{2}<\cdots <i_{r}\leq n$. Let
\[ I=<u>_{\de} = \prod_{q=1}^{r}\prod_{j=0}^{s}(\me_{q}^{[d_{j}]})^{\alpha_{qj}},\] where $\alpha_{q} = \sum_{j=0}^{s}\alpha_{qj}d_{j}$. Suppose $i_{r}=n$. Let $1\leq a \leq r$ be an integer and
\[ P_{a}(I):=\{ (\lambda,t)\in \mathbb N^{a}\times \mathbb N^{a}|
\;1\leq \lambda_{1}<\cdots < \lambda_{a}=r,
t_{a}>\cdots >t_{1}, \alpha_{\lambda_{\nu} t_{\nu}}\neq 0,\; for\; 1\leq \nu \leq a \}.  \]
Let $J = \sum_{a=1}^{r}\sum_{(\lambda,t)\in P_{a}(I)}J_{(\lambda,t)}$, where $J_{(\lambda,t)}$ is the ideal
\[ \prod_{e=1}^{a}(x_{i_{\lambda_{e}}} \cdots x_{i_{\lambda_{e-1}}+1})^{d_{t_{e}}-1}
\prod_{\nu=1}^{a}\me_{\lambda_{\nu}}^{[d_{t_{\nu}+1}]}\prod_{j>t_{\nu}} (\me_{\lambda_{\nu}}^{[d_{j}]})^{\alpha_{\lambda_{\nu}j}} (\me_{\lambda_{\nu}}^{[d_{t_{\nu}}]})^{\alpha_{\lambda_{\nu},t_{\nu}}-1} \prod_{q=\lambda_{\nu-1}+1}^{\lambda_{\nu}-1} \prod_{j\geq t_{\nu}} (m_{q}^{[d_{j}]})^{\alpha_{qj}}, \]
where we denote $\me^{[d_{t_{a+1}}]} = S$. Then $Soc(S/I) = (J+I)/I$.
\end{teor}

\begin{proof}
The proof will be given by induction on $r$, the case $r=1$ being done in Lemma $2.1$. Suppose that $r>1$. For $1\leq q \leq r$, let: $I_{q}=\prod_{e=1}^{q}\prod_{j=0}^{s}(\me_{e}^{[d_{j}]})^{\alpha_{ej}}$ and $S_{q}=k[x_{1},x_{2},\ldots,x_{i_{q}}]$
For $t$ with $\alpha_{rt}\neq 0$, denote:
\[ I^{(t)} = \me_{r-1}^{[d_{t}]}\prod_{j<t}(\me_{r-1}^{[d_{j}]})^{\alpha_{r-1,j}}I_{r-2}.\]
Let $J^{(t)}$ be an ideal in $S_{r-1}$ such that $Soc(S_{r-1}/I^{(t)}) = (J^{(t)}+I^{(t)})/I^{(t)}$. The induction step is given in the following lemma:
\begin{lema}
Suppose $i_{r}=n$ and let
\[ J = \sum_{t=0,\alpha_{rt}\neq 0} (x_{n}\cdots x_{i_{r-1}+1})^{d_{t}-1} \prod_{j>t}(\me_{r}^{[d_{j}]})^{\alpha_{r}j}
 \prod_{j\geq t}(\me_{r-1}^{[d_{j}]})^{\alpha_{r-1,j}}(\me_{r}^{[d_{t}]})^{\alpha_{rt}-1}J^{(t)}. \]
 Then $Soc(S/I) = (J+I)/I$.
\end{lema}
\begin{proof}
Let $w \in S\setminus I$ be a monomial such that $\me_{r}w \subset I$. As in the proof of lemma $2.1$, we choose for each $1\leq \rho \leq n$, $e_{\rho}= max\{e:\; x_{\rho}^{d_{e}-1}|w\}$. Renumbering variables $\{x_{n}, \ldots, x_{i_{r-1}+1} \}$ (it does not affect $I$, $J$ and $I^{(t)}$), we may suppose $e_{n}\leq e_{n-1}\leq \cdots \leq e_{i_{r-1}+1}$. Set $t=e_{n}$. We claim that $\alpha_{rt}\neq 0$. Indeed, if $\alpha_{rt}=0$ then from $x_{n}w\in I$ we get $x_{n}w/x_{n}^{d_{t}-1} \in \widetilde{I} = \prod_{j>t}(\me_{r}^{[d_{j}]})^{\alpha_{rj}}I_{r-1}$ because $x_{n}^{d_{t}-1} \in \prod_{j<t}(\me_{r}^{[d_{j}]})^{\alpha_{rj}}$. Since $t=e_{n}$ is maximal chosen, we get $w/x_{n}^{d_{t}-1} \in \widetilde{I}$ and so $w\in I$ a contradiction.

Reduction to the case that $x_{n}^{d_{t}}$ does not divide $w$. Suppose that
$w =x_{n}^{d_{t}}\widetilde{w}$ and set
\[ \widetilde{I} = (\me_{r}^{[d_{t}]})^{\alpha_{rt}-1} \prod_{\epsilon \leq 0,\epsilon\neq t}(\me_{r}^{[d_{\epsilon}]})^{\alpha_{r\epsilon}}I_{r-1}.\]
We see that $\me w\in I \Leftrightarrow \me \widetilde{w}\in \widetilde{I}$. Replacing $w$ and $I$ with $\widetilde{w}$ and $\widetilde{I}$, we reduce our problem to a new $\widetilde{t}<t$. The above argument implies that $\widetilde{\alpha}_{r\widetilde{t}} \neq 0$, where $\widetilde{\alpha}$ is the 'new' $\alpha$ of $\widetilde{I}$.

Reduction to the case when $\alpha_{rj} = \alpha_{r-1,j} = 0$ for $j>t$, $\alpha_{rt}=1$ and $\alpha_{r-1,t}=0$. From $x_{n}w\in I$, we see that there exists $\rho<n$ such that $x_{\rho}^{d_{j}}|w$ for $j>t$ if $\alpha_{rj}\neq 0$, or $j=t$ if $\alpha_{rt}>1$. Choose such maximal possible $\rho$. Set $w' = w/x_{\rho}^{d_{j}}$,
\[ I' = (\me_{r}^{[d_{j}]})^{\alpha_{rj}-1} \prod_{\epsilon\geq 0, \epsilon\neq j} (\me_{r}^{[d_{\epsilon}]})^{\alpha_{r\epsilon}}I_{r-1}. \]
We see that $\me w \subset I \Leftrightarrow \me w' \subset I'$, because from $x_{n}w\in I$, we get $x_{n}w'\in I'$ from the maximality of  $\rho$.

Let $\alpha'_{rj}=\alpha_{rj}-1$ and $\alpha'_{q\epsilon} = \alpha_{q\epsilon}$ for $(q,\epsilon) \neq (r,j)$. $\alpha'$ is the 'new' $\alpha$ for $I'$. If we show that
\[ w'\in J' = \sum_{e \geq 0, \alpha'_{re}\neq 0} (x_{n}\cdots x_{i_{r-1}+1})^{d_{e}-1} \prod_{\epsilon >e} (\me_{r}^{[d_{\epsilon}]})^{\alpha_{r}\epsilon}
 \prod_{j\geq e}(\me_{r-1}^{[d_{j}]})^{\alpha_{r-1,j}}(\me_{r}^{[d_{\epsilon}]})^{\alpha_{re}-1}J^{(t)}, \]
then $w = x_{\rho}^{d_{j}}w' \in \me_{r}^{[d_{j}]}J'\subset J$.
Using this procedure, by recurrence we arrive to the case $\alpha_{rj}=0$ for $j>t$ and $\alpha_{rt}=1$. Again from $x_{n}w \in I$, we note that there exists $\rho < i_{r-1}$ such that $x_{\rho}^{d_{j}}|w$ for $j\geq t$ with $\alpha_{r-1,j}\neq 0$. Choose such maximal possible $\rho$ and note that $\me w \subset I$ if and only if $\me w'' \in I''$ for $w'' = w/ x_{\rho}^{d_{j}}$, where \[ I'' = (\me_{r-1}^{d_{j}})^{\alpha_{r-1,j-1}} \prod_{\epsilon\geq 0,\;\epsilon\neq j} (\me_{r-1}^{[d_{\epsilon}]})^{\alpha_{r-1, \epsilon}} \prod_{\epsilon \geq 0} (\me_{r}^{[d_{\epsilon}]})^{\alpha_{r\epsilon}}I_{r-2}.\]
As above, we reduce our problem to $I''$ and the $\alpha''$, which is the new $\alpha$ of $I''$, is given by
$\alpha''_{r-1,j} = \alpha_{r-1,j-1}$, $\alpha''_{q\epsilon} = \alpha_{q\epsilon}$ for $(q,\epsilon) \neq (r-1,j)$.
Using this procedure, by recurrence we end our reduction.

Case $\alpha_{rj} = \alpha_{r-1j} = 0$ for $j>t$, $\alpha_{rt}=1$ , $\alpha_{r-1t}=0$ and $x_{n}^{d_{t}}$ does not divide $w$. Let express $w=(x_{n}\cdots x_{i_{r-1}+1})^{d_{t}-1}y$. We will show that $y$ does not depend on $\{x_{n},\ldots,x_{i_{r-1}+1}\}$. Indeed, if $n=i_{r-1}+1$ then there is nothing to show since $x_{n}^{d_{t}}$ does not divide $w$. Suppose that $n>i_{r-1}+1$, then from $x_{n}w\in I$ we get $y\in I_{r-1}$ because $x_{n-1}^{d_{t}-1}\in \prod_{j<t} (\me_{r}^{d_{j}})^{\alpha_{rj}}$ and the variables $x_{n}, \ldots, x_{i_{r-1}+1}$ are regular on $S/I_{r-1}S$. If $y=x_{\eta}y'$ for $\eta>i_{r-1}$, then as above $y' \in I_{r-1}$. Thus $w\in x_{\eta}^{d_{t}}x_{\rho}^{d_{t}-1}y' \subset I$ for any $\rho \neq \eta, i_{r-1} < \rho \leq n$, a contradiction.

Note that $\me_{r}w\in I \Rightarrow \me_{r-1}y \in I^{(t)}$ and so $w\in (x_{n}\cdots x_{i_{r-1}+1})^{d_{t}-1}J^{(t)}$. Since $\alpha_{rj} = \alpha_{r-1j} = 0$ for $j>t$ and $\alpha_{rt}=1$ and $\alpha_{r-1,t} = 0$, we get $w\in J$. Conversely, if $y\in J^{(t)}$, then it is clear that $w\in J$.
\end{proof}

We see by the above lemma that:
\[(*)\;J = \sum_{e\geq 0,\alpha_{re}\neq 0} (x_{n} \cdots x_{i_{r-1}+1})^{d_{e}-1} \prod_{j>e} (\me_{r-1}^{[d_{j}]})^{\alpha_{rj}} \prod_{j\geq e} (\me_{r-1}^{[d_{j}]})^{\alpha_{r-1,j}} (\me_{r}^{[d_{e}]})^{\alpha_{re}-1}J^{(e)}.\]
Since $\lambda_{a}=r-1$, by the induction hypothesis applied to $I^{(e)}$ we get:
\[ J^{(e)} = \sum_{a=1}^{r-1} [\sum_{(\lambda,t)\in P_{a}(I^{(e)}),t_{a}=e} \prod_{s=1}^{a} (x_{i_{\lambda_{s}}} \cdots x_{i_{\lambda_{s-1}+1}} )^{d_{t_{s}}-1} \cdot J'_{(\lambda,t)} + \]\[ +
\sum_{(\lambda,t)\in P_{a}(I^{(e)}),t_{a}<e} \prod_{s=1}^{a} (x_{i_{\lambda_{s}}} \cdots x_{i_{\lambda_{s-1}+1}} )^{d_{t_{s}}-1} \cdot J''_{(\lambda,t)} ],\;where\; \]
\[ J'_{(\lambda,t)} = \prod_{q=\lambda_{a-1}+1}^{\lambda_{a}-1}\prod_{j\geq e} (\me_{q}^{[d_{j}]})^{\alpha_{qj}}\widetilde{J}_{(\lambda,t)} \;and\;\]
\[ J''_{(\lambda,t)} = \me_{r-1}^{[d_{e}]} \prod_{j>t_{a}}^{e-1}(\me_{\lambda_{a}}^{[d_{j}]})^{\alpha_{\lambda_{a},j}} (\me_{\lambda_{a}}^{[d_{t_{a}}]})^{\alpha_{\lambda_{a},t_{a}-1}}\cdot \prod_{q=\lambda_{a-1}+1}^{\lambda_{a}-1} \prod_{j\geq t_{a}} (\me_{q}^{[d_{j}]})^{\alpha_{qj}}\widetilde{J}_{(\lambda,t)} ,\;and\]
\[ \widetilde{J}_{(\lambda,t)} = \prod_{\nu=1}^{a-1}\me_{\lambda_{\nu}}^{[d_{t_{\nu+1}}]}
\prod_{j>t_{\nu}}(\me_{\lambda_{\nu}}^{[d_{j}]})^{\alpha_{\lambda_{\nu},j}} (\me_{\lambda_{\nu}}^{[d_{t_{\nu}}]})^{\alpha_{\lambda_{\nu},t_{\nu}-1}}\cdot
 \prod_{q=\lambda_{\nu-1}+1}^{\lambda_{\nu}-1}\prod_{j\geq t_{\nu}} (\me_{q}^{[d_{j}]})^{\alpha_{qj}}.\]

If $t_{a}=e$, set $\lambda'_{\nu} = \lambda_{\nu}$ for $\nu<a$, $\lambda'_{a}=r$ and see that $(\lambda',t) \in P_{a}(I)$. If $t_{a}<e$, then put $\lambda''_{\nu} = \lambda_{\nu}$ for $\nu\leq a$, $\lambda''_{a+1}=r$, $t''_{\nu} = t_{\nu}$ for $\nu\leq a$ and $t''_{a+1}=e$ and then $(\lambda'',t) \in P_{a+1}(I)$. Substituting $J^{(e)}$ in $(*)$, we get the following expression for $J$:
\[\sum_{a=1}^{r-1} \sum_{(\lambda',t)\in P_{a}(I)} \prod_{\nu=1}^{a} (x_{i_{\lambda'_{\nu}}}\cdots x_{i_{\lambda'_{\nu-1}+1}})^{d_{t_{\nu}}-1} \cdot [ \prod_{j>e} (\me_{\lambda'_{a}}^{[d_{j}]})^{\alpha_{\lambda'_{a}j}} (\me_{\lambda'_{a}}^{[d_{e}]})^{\alpha_{\lambda'_{a}e}-1}
 \cdot \prod_{q=\lambda'_{a-1}+1}^{\lambda'_{a}-1}\prod_{j\geq e} (\me_{q}^{[d_{j}]})^{\alpha_{qj}}]\cdot \] \[ \cdot \widetilde{J}_{(\lambda,t)} +
\sum_{a=1}^{r-1} \sum_{(\lambda'',t'')\in P_{a+1}(I)} \prod_{\nu=1}^{a+1} (x_{i_{\lambda''_{\nu}}}\cdots x_{i_{\lambda''_{\nu-1}+1}})^{d_{t_{\nu}}-1} \cdot  [\prod_{j>e}(\me_{\lambda''_{a+1}}^{[d_{j}]})^{\alpha_{\lambda''_{a+1},j}}
(\me_{\lambda''_{a+1}}^{[d_{t''_{a}+1}]})^{\alpha_{\lambda''_{a+1}t''_{a+1}}-1}] \]\[
[ \me_{\lambda''_{a}}^{[d_{t''_{a+1}}]} \prod_{j\geq t''_{a}} (\me_{\lambda''_{a}}^{[d_{j}]})^{\alpha_{\lambda''_{a+1},j}}
(\me_{\lambda''_{a}}^{[d_{t''_{a}}]})^{\alpha_{\lambda''_{a}t''_{a}}-1}
\prod_{q=\lambda''_{a-1}+1}^{\lambda''_{a}-1}\prod_{j\geq t''_{a}} (\me_{q}^{[d_{j}]})^{\alpha_{qj}}] \cdot \widetilde{J}_{(\lambda,t)}. \]
Since all the pairs of $P_{b}(I)$ have the form $(\lambda',t)$ or $(\lambda'',t'')$ for a pair $(\lambda,t)\in P_{b}(I)$ or $(\lambda,t)\in P_{b-1}(I)$ respectively, it is not hard to see that the expression above is the formula of $J$ as stated.
\end{proof}

\noindent
Let $s_{q}=max\{j|\alpha_{qj}\neq 0\}$, $d_{qt} = \sum_{e=1}^{q}\sum_{j\geq t}^{s_{q}}\alpha_{ej}d_{j}$,
$D_{q}=d_{q,s_{q}} + (i_{q}-1)(d_{s_{q}}-1)$ for $1\leq q\leq r$.

\begin{cor}
With the notation and hypothesis of above theorem, for $(\lambda,t)\in P_{a}(I)$ let:
\[d_{(\lambda,t)} = \sum_{\nu=1}^{a} \sum_{q=\lambda_{\nu-1}+1}^{\lambda_{\nu}} \sum_{j\geq t_{\nu}} \alpha_{qj}d_{j}.\; Then:\]
\begin{enumerate}
    \item $Soc(I_{r-1}S/I) = Soc(S/I)$.
    \item $((J+I)/I)_{e}\neq 0$, if and only if
          $ e = d_{(\lambda,t)} + \sum_{\nu=1}^{a}(i_{\lambda_{\nu}} - i_{\lambda_{\nu-1}})(d_{t_{\nu}}-1) - d_{t_{1}}$,
          for some $1\leq a\leq r$ and $(\lambda,t)\in P_{a}(I)$.
    \item $c=max\{e|((J+I)/I)_{e}\neq 0  \} = d_{r,s_{r}} + (n-1)(d_{s_{r}} -1) - 1$.
\end{enumerate}
\end{cor}

\begin{proof}
1.Note that $J_{(\lambda,t)}$ is contained in
\[ \prod_{q=1,q\notin \{\lambda_{1},\ldots,\lambda_{q}\}}^{r} (\me_{q}^{[d_{j}]})^{\alpha_{qj}}
\prod_{\nu=1}^{a}[\prod_{j\neq t_{\nu}} (\me_{\lambda_{\nu}}^{d_{j}})^{\alpha_{\lambda_{\nu},j}}
(\me_{\lambda_{\nu}}^{t_{\nu}})^{\alpha_{\lambda_{\nu},t_{\nu}-1}}]\prod_{\epsilon=1}^{a-1}\me_{\lambda_{\epsilon+1}}^{d_{t_{\epsilon}}+1}.\]
Since $\me_{\lambda_{\epsilon}}^{d_{t_{\epsilon}+1}} \subset \me_{\lambda_{\epsilon}}^{d_{t_{\epsilon}}}$ for $t_{\epsilon+1}>t_{\epsilon}$ and $\lambda_{a}=r$ if follows that
\[ J \subset \prod_{j\neq t_{a}} (\me_{r}^{[d_{j}]})^{\alpha_{rj}} (\me_{r}^{[d_{t_{a}}]})^{\alpha_{rt_{a}}-1}I_{r-1},\]
as desired.

2.If $((J+I)/I)_{e}\neq 0$ then there exists a monomial $u\in J\setminus I$ of degree $e$. But $u\in J$, implies that there exists $a\in \{1,\ldots,r\}$ and $(\lambda,t)\in P_{a}(I)$ such that $u\in J_{(\lambda,t)}$. Thus the degree of $u$ is $e = d_{(\lambda,t)} + \sum_{\nu=1}^{a}(i_{\lambda_{\nu}} - i_{\lambda_{\nu-1}})(d_{t_{\nu}}-1) - d_{t_{1}}$, as required.

Conversely, let $e = d_{(\lambda,t)}+ \sum_{\nu=1}^{a}(i_{\lambda_{\nu}} - i_{\lambda_{\nu-1}})(d_{t_{\nu}}-1) - d_{t_{1}}$ for some $a\in \{1,\ldots,r\}$ and $(\lambda,t)\in P_{a}(I)$. We show that the monomial
\[w = \prod_{\nu=1}^{a} (x_{i_{\lambda_{\nu}}}\cdots x_{i_{\lambda_{\nu-1}+1}})^{d_{t_{\nu}}-1}\cdot x_{1}^{d_{(\lambda,t)}-d_{t_{1}}} \in J\setminus I. \]
Obvious $w\in J$. Let us assume that $w\notin I$. Then $w/x_{i_{\lambda_{a}}}^{d_{t_{a}}-1} \in \prod_{j\geq t_{a}} (\me_{\lambda_{a}}^{[d_{j}]})^{\alpha_{\lambda_{a}j}}I_{\lambda_{a}-1}$ because $x_{i_{\lambda_{a}}}^{d_{t_{a}}-1} \in \prod_{j<t_{a}}(\me_{\lambda_{a}}^{[d_{j}]})^{\alpha_{\lambda_{a}j}}$ and $x_{i_{\lambda_{a}}} \notin \me_{j}$ for $j<\lambda_{a}$. Inductively we get that:
\[ w/(x_{i_{\lambda_{a}}} \cdots x_{i_{\lambda_{a-1}+1}})^{d_{t_{a}}-1} \in \prod_{q=\lambda_{a-1}+1}^{\lambda_{a}}  \prod_{j\geq t_{a}} (\me_{\lambda_{a}}^{[d_{j}]})^{\alpha_{qj}}I_{\lambda_{a}-1}. \]
Following the same reduction and using that $t_{a}>\cdots>t_{1}$ we obtain that:
\[ x_{1}^{d_{(\lambda,t)}-d_{t_{1}}} \in \prod_{\nu=1}^{a} \prod_{q=\lambda_{\nu-1}+1}^{\lambda_{\nu}}  \prod_{j\geq t_{\nu}} (\me_{\lambda_{a}}^{[d_{j}]})^{\alpha_{qj}}. \]
So $d_{(\lambda,t)}-d_{t_{1}} \geq d_{(\lambda,t)}$, a contradiction.

3.Note that $c=d_{(\lambda',t')}$ for $(\lambda',t')\in P_{1}(I)$ with $\lambda' = \lambda_{1} = r$ and $t'=t_{1}=s_{r}$. We have to show that:
\[ c = d_{r,s_{r}} + (n-1)(d_{s_{r}} -1) - 1 \leq d_{(\lambda,t)} + \sum_{\nu=1}^{a}(i_{\lambda_{\nu}} - i_{\lambda_{\nu-1}})(d_{t_{\nu}}-1) - d_{t_{1}}, \]
for any $1\leq a\leq r$ and $(\lambda,t)\in P_{a}(I)$. Since $d_{s_{r}}-1 \leq (d_{t_{\nu}}-1) + \sum_{j\leq t_{\nu}}^{s_{r}-1}\alpha_{qj}d_{j}$ for all $q$ with $i_{\nu-1}<q \leq i_{\nu}$, we see that:
\[ d_{r,s_{r}} + (n-1)(d_{s_{r}} -1)-1 \geq d_{(\lambda,t)} + \sum_{\nu=2}^{a}(i_{\lambda_{\nu}} - i_{\lambda_{\nu-1}})(d_{t_{\nu}}-1) + (i_{\lambda_{1}}-1)(d_{t_{1}}-1) -1. \]
On the other hand, $(i_{\lambda_{1}}- i_{\lambda_{0}})(d_{t_{1}}-1) = (i_{\lambda_{1}}- 1)(d_{t_{1}}-1) + d_{t_{1}}-1$,
and replacing that in the above relation we obtained what is required.
\end{proof}

\begin{exm}
Let $\de: 1|2|4|12$.
\begin{enumerate}
    \item Let $u=x_{3}^{21}$. We have $\alpha_{0}=1$, $\alpha_{1}=0$,  $\alpha_{2}=2$ and $\alpha_{3}=1$ so:
    \[ I = <u>_{\de} = (x_{1},x_{2},x_{3})(x_{1}^{4},x_{2}^{4},x_{3}^{4})^{2}(x_{1}^{12},x_{2}^{12},x_{3}^{12}).\]
    Let $J= \sum_{t=0,\alpha_{t}>0} J_{t}$, where $J_{t} = (x_{1}x_{2}x_{3})^{d_{t}-1}(x_{1}^{d_{t}},
                x_{2}^{d_{t}}, x_{3}^{d_{t}})^{\alpha_{t}-1}\prod_{j>t}(x_{1}^{d_{j}},
                x_{2}^{d_{j}}, x_{3}^{d_{j}})^{\alpha_{j}}$.

    $J_{0} = (x_{1}x_{2}x_{3})^{1-1}\cdot (x_{1},x_{2},x_{3})^{1-1} \cdot \prod_{j>t}(x_{1}^{d_{j}},
                x_{2}^{d_{j}}, x_{3}^{d_{j}})^{\alpha_{j}} = (x_{1}^{4},x_{2}^{4},x_{3}^{4})^{2}(x_{1}^{12},x_{2}^{12},x_{3}^{12}).$

    $J_{2} =(x_{1}x_{2}x_{3})^{4-1}(x_{1}^{4},x_{2}^{4},x_{3}^{4})^{2-1}(x_{1}^{12},x_{2}^{12},x_{3}^{12}) = (x_{1}x_{2}x_{3})^{3}(x_{1}^{4},x_{2}^{4},x_{3}^{4})(x_{1}^{4},x_{2}^{4},x_{3}^{12})$ and

    $ J_{3}=(x_{1}x_{2}x_{3})^{12-1} = (x_{1}x_{2}x_{3})^{11}.$ From $2.1$ , $Soc(S/I) = (J+I)/I$.

    \item Let $u=x_{2}^{9}x_{3}^{16}$. We have $r=2$, $i_{1}=2$ and $i_{2}=3$. Also $\alpha_{10}=1$, $\alpha_{12}=2$,
          $\alpha_{22}=1$, $\alpha_{23}=1$ and the other components of $\alpha$ are zero. Then
    \[ I = <u>_{\de} = <x_{2}^{9}>_{\de}<x_{3}^{16}>_{\de} = (x_{1},x_{2})(x_{1}^{4},x_{2}^{4})^{2}
                    (x_{1}^{4},x_{2}^{4},x_{3}^{4})(x_{1}^{12},x_{2}^{12},x_{3}^{12}).\]
    We have two possible partitions: (a) $(2)$ and (b) $(1<2)$.

    (a)$\lambda=\lambda_{1}=2$, $t=t_{1}$ such that $\alpha_{2t}\neq 0$. We have two possible $t$: $t=2$ or $t=3$.

    (i)For $t=2$ we obtain (according to the Theorem $2.4$) the following part of the socle:
    \[ J_{(2,2)} = (x_{1}x_{2}x_{3})^{3} (x_{1}^{12}, x_{2}^{12}, x_{3}^{12}) (x_{1}^{4},x_{2}^{4})^{4} \]

  (ii)For $t=3$ we obtain:
    \[ J_{(2,3)} = (x_{1}x_{2}x_{3})^{11}  \]

  (b)$1=\lambda_{1}<\lambda_{2}=2$, $t=(t_{1},t_{2})$ such that $\alpha_{\lambda_{e},t_{e}}\neq 0$ for $1\leq e\leq 2$
  and $t_{1}<t_{2}$. According to our expressions for $\alpha_{i}$ we have three possible cases: $t_{1}=0,t_{2}=2$ or
  $t_{1}=0,t_{2}=3$ or $t_{1}=2,t_{2}=3$.

  (i)For $t_{1}=0$ and $t_{2}=2$ we obtain:
  \[ J_{(1,2),(0,2)} = x_{3}^{3}(x_{1}^{4},x_{2}^{4})(x_{1}^{4},x_{2}^{4})^{2}(x_{1}^{12},x_{2}^{12},x_{3}^{12}). \]

  (ii)For $t_{1}=0$ and $t_{2}=3$ we obtain:
  \[ J_{(1,2),(0,3)} = x_{3}^{11}(x_{1}^{12},x_{2}^{12})(x_{1}^{4},x_{2}^{4})^{2} \]

  (iii)For $t_{1}=2$ and $t_{2}=3$ we obtain:
  \[ J_{(1,2),(2,3)} =  x_{1}^{3}x_{2}^{3}x_{3}^{11}(x_{1}^{12},x_{2}^{12})(x_{1}^{4},x_{2}^{4})\]

  From $2.4$ it follows that if $J = J_{(2,2)} + J_{(2,3)} + J_{(1,2),(0,2)} + J_{(1,2),(0,3)} + J_{(1,2),(2,3)}$ then
  $Soc(S/I)=(I+J)/J$.
\end{enumerate}
\end{exm}

\vspace{9mm}
\section{A generalization of Pardue's formula.}

In this section, we give a generalization of a theorem proved by Aramova-Herzog \cite{ah} and Herzog-Popescu \cite{hp} which is known as "Pardue's formula".

Let $1\leq i_{1}<i_{2}<\cdots<i_{r}=n$ and let $\alpha_{1},\ldots,\alpha_{r}$ some positive integers.
Let $u=\prod_{i=1}^{r}x_{i_{q}}^{\alpha_{q}} \in S=K[x_{1},\ldots,x_{n}]$. Our goal is to give a formula for the regularity of the ideal
\[ I = <u>_{\de} = \prod_{r=1}^{q}\prod_{j=0}^{s}(\me_{q}^{[d_{j}]})^{\alpha_{qj}}, \]
where $\alpha_{q}=\sum_{j=0}^{s}\alpha_{qj}d_{j}$. If $i_{1}=1$, it follows that $I=x_{1}^{\alpha_{1}}I'$, where
$I' = \prod_{r=2}^{q}\prod_{j=0}^{s}(\me_{q}^{[d_{j}]})^{\alpha_{qj}}$, and therefore $reg(I) = \alpha_{1} + reg(I')$. Thus, we may assume $i_{1}\geq 2$.

If $N$ is a graded $S$-module of finite length, we denote $s(N)=max\{i|N_{i}\neq 0\}$. Let $s_{q}=max\{j|\alpha_{qj}\neq 0\}$ and $d_{qt}= \sum_{e=1}^{q}\sum_{j\geq t}^{s_{e}}\alpha_{ej}d_{j}$.
Let $D_{q} = d_{qs_{q}} + (i_{q}-1)(d_{s_{q}}-1)$, for $1\leq q \leq r$. With this notations we have:

\begin{teor}
$reg(I) = max_{1\leq q \leq r} D_{q}$. In particular, if $I = <x_{n}^{\alpha}>_{\de}$ and $\alpha=\sum_{t=0}^{s} \alpha_{t}d_{t}$ with $\alpha_{s}\neq 0$ then $reg(I) = \alpha_{s}d_{s} + (n-1)(d_{s}-1)$.
\end{teor}

\begin{proof}
Let $I_{\ell}=\prod_{q=1}^{r-\ell}\prod_{j=0}^{s}(\me_{q}^{[d_{j}]})^{\alpha_{qj}}$, for $0 \leq \ell \leq r$.
Then $I=I_{0}\subset I_{1}\subset \cdots \subset I_{r} = S$ is the sequential chain of ideals of $I$, i.e.
$I_{\ell+1} = (I_{\ell}:x_{n_{\ell}}^{\infty})$, where $n_{\ell} = i_{r-\ell}$. Moreover, from the Remark $1.13$, we see that this chain is in fact the chain from the definition of a sequentially Cohen-Macauly module for $S/I$.
Let $S_{\ell}=k[x_{1},\ldots,x_{n_{\ell}}]$ and $m_{\ell} = (x_{1},\ldots,x_{n_{\ell}})$.

The corollary $2.6$ implies that $c_{e}=D_{e}-1$ is the maximal degree for a nonzero element of $Soc(S_{\ell}/J_{\ell})$. \cite[Corollary 2.7]{hpv} implies $reg(I) = max\{s(I_{\ell}S_{\ell}^{sat}/I_{\ell}S_{\ell})\;|\; \ell=0,\ldots,r-1)\} + 1$. Also, from the corollary $2.6$, we get
\[ Soc(S_{\ell}/I_{\ell}S_{\ell}) = Soc(I_{\ell+1}S_{\ell}/I_{\ell}S_{\ell}) = (I_{\ell+1}:m_{\ell})S_{\ell}/I_{\ell}S_{\ell} = I_{\ell}S_{\ell}^{sat}/I_{\ell}S_{\ell}, \]
which complete the proof.
\end{proof}

\begin{cor}
$reg(I) \leq n\cdot deg(u) = n \cdot deg(I)$, where $deg(I)=max\{deg(w)|w\in G(I)\}$.
\end{cor}

\begin{cor}
$S/I$ has at most $r$-corners among $(i_{q},D_{q}-1)$ for $1\leq q\leq r$. If $i_{1}=1$ we replace $(i_{1},D_{1}-1)$ with $(1,\alpha_{1})$. The corresponding extremal Betti numbers are $\beta_{i_{q},D_{q}+i_{q}-1}$.
\end{cor}

\begin{proof}
By Theorem $1.14$ combined with the proof of Theorem $3.1$, $S/I$ has at most $r$-corners among
$(n_{\ell},s(I_{\ell+1}S_{\ell}/I_{\ell}S_{\ell}))$ and is enough to apply Corollary $2.6$.
\end{proof}

\begin{exm}
Let $\de: 1|2|4|12$.
\begin{enumerate}
    \item Let $u=x_{3}^{21} \in k[x_{1},x_{2},x_{3}]$. We have $21=1\cdot 1 + 0\cdot 2 + 2 \cdot 4 + 1 \cdot 12$. From
    $3.1$, we get: \[ reg(<u>_{\de}) = 1\cdot 12 + (3-1)\cdot (12 - 1) = 34.\]
    \item Let $u=x_{1}^{2}x_{2}^{16}x_{3}^{9}$. Then $reg(<u>_{\de}) = 2 + reg(<u'>_{\de})$, where $u' = u/x_{1}^{2}$.
          We compute $reg(<u'>_{\de})$. With the notations above, we have $i_{1}=2$, $i_{2}=3$, $r=2$, $\alpha_{1}=16$ and
          $\alpha_{2}=9$. We have $\alpha_{1} = 1\cdot 4 + 1\cdot 12$ and $\alpha_{2} =1\cdot 1 + 2\cdot 4$, thus
          $s_{1}=3$ and $s_{2}=2$. $D_{1} = d_{13} + (2-1)(d_{3}-1) = 12 + 11 = 23$ and
          $D_{2} = d_{22} + (3-1)(d_{2}-1) = 24 + 6 = 30$. In conclusion, $reg(<u>_{\de}) = 2 + max\{23,30\} = 32$.
\end{enumerate}
\end{exm}

In the following, we show that if $I$ is a principal
$\de$-fixed ideal generated by the power of a variable, then $I_{\geq e}$ is stable for any $e\geq reg(I)$.

\begin{lema}
Let $I=<x_{n}^{\alpha}>_{\de}$ and $\alpha=\sum_{t=0}^{s} \alpha_{t}d_{t}$ with $\alpha_{s}\neq 0$. If $e\geq reg(I)+1$ then for every monomial $v\in I_{\geq e}$ there exists $w\in G(I)$ and a monomial $y\in S$ such that $v=w\cdot y$ and $m(v)=m(y)$.
\end{lema}

\begin{proof}
We may assume $e=reg(I)+1$ and $v\in I_{e}$. Then $v=w'\cdot y'$ for some $w'\in G(I)$ and a monomial $y'\in S$. Suppose
$w' = \prod_{t=0}^{s}\prod_{j=1}^{n} x_{j}^{\lambda_{tj}\cdot d_{t}}$, where $0\leq \lambda_{tj}$ and $\sum_{j=1}^{n}\lambda_{tj} = \alpha_{t}$. Suppose $n = m(v) = m(w') > m(y')$. Then $y'=x_{1}^{\beta_{1}}\cdots x_{n-1}^{\beta_{n-1}}$. Let $m = min\{t| \lambda_{tn}\neq 0 \}$.

We claim that there exists some $1\leq i \leq n$ such that $d_{m} - \sum_{t=0}^{m-1}\lambda_{ti}d_{t} \leq \beta_{i}$.
Otherwise, it follows that $d_{m} - \sum_{t=0}^{m-1}\lambda_{ti}d_{t} \geq \beta_{i} + 1$ for any $i=1,\ldots,n-1$. So,
\[ (n-1)d_{m} - \sum_{i=1}^{n-1}\sum_{t=0}^{m-1}\lambda_{ti}d_{t} \geq \beta_{1} + \cdots + \beta_{n-1} + n - 1 = reg(I) + 1 - \alpha + n-1 \Leftrightarrow \]
\[ (n-1)(d_{m}-1) - \sum_{t=0}^{m-1}\alpha_{t}d_{t} \geq (n-1)(d_{s}-1) - \sum_{t=0}^{s-1}\alpha_{t}d_{t} + 1 \Leftrightarrow  \]
\[ \sum_{t=m}^{s-1}\alpha_{t}d_{t}\geq (n-1)(d_{s}-d_{m}) + 1, \]
because $reg(I)=\alpha_{s}d_{s} + (n-1)(d_{s}-1)$ from Theorem $3.1$. But on the other hand, $d_{s}-d_{m} = \sum_{t=m}^{s-1} (d_{t+1}/d_{t}-1)d_{t}  \geq \sum_{t=m}^{s-1}\alpha_{t}d_{t}$ and this contradict the above inequality.

Thus, we may choose $i<n$ such that $\gamma = d_{m} - \sum_{t=0}^{m-1}\lambda_{ti}d_{t} \leq \beta_{i}$. Therefore, we can write: $v = w'\cdot y' = w\cdot y$, where $w=w'\cdot x_{i}^{\gamma}/x_{n}^{\gamma}$ and $y = w'\cdot x_{n}^{\gamma}/x_{i}^{\gamma}$. It is easy to see that $w \in G(I)$ and $m(v)=m(y)=n$.
\end{proof}

\begin{cor}
If $I=<x_{n}^{\alpha}>_{\de}$ and $e\geq reg(I)$ then $I_{\geq e}$ is stable.
\end{cor}

\begin{proof}
Let $v\in I_{\geq e}$. Let $i<m(v)$. Since $x_{i}\cdot v \in I_{\geq e+1}$ it follows from the above lemma that
$x_{i}v=w\cdot y$ for some $w\in G(I)$ and $y \in S$ such that $m(x_{i}v) = m(y)$. But $m(v)=m(x_{i}v)$ and thus
$x_{i}v/x_{m(v)} = w \cdot y/x_{m(v)} \in I$.
\end{proof}

The converse is also true. Indeed we have the following more general result of Eisenbud-Reeves-Totaro:

\begin{prop}\cite[Proposition 12]{ert}
Let $I$ be a monomial ideal with $deg(I)=d$ and let $e\geq d$ such that $I_{\geq e}$ is stable. Then $reg(I)\leq e$.
\end{prop}

\begin{obs}
$3.6$ gives another proof for the "$\leq$" inequality of the generalised Pardue's formula in the case when $I = <x_{n}^{\alpha}>_{\de}$. Indeed, considering $e=\alpha_{s}d_{s}+(n-1)(d_{s}-1)$ from $3.6$ it follows that $I_{\geq e}$ is stable and thus $3.7$ implies $reg(I)\leq e$.
\end{obs}

\begin{cor}
If $I=<x_{n}^{\alpha}>_{\de}$ then $reg(I)=min\{e|\;I_{\geq e}$ is stable $\}$.
\end{cor}

\vspace{2mm} \noindent {\footnotesize
\begin{minipage}[b]{10cm}
 Mircea Cimpoea\c s, Junior Researcher\\
 Institute of Mathematics of the Romanian Academy\\
 Bucharest, Romania\\
 E-mail: mircea.cimpoeas@imar.ro

\end{minipage}}


\begin{thebibliography}{breitestes Label}
  \bibitem[1]{ah}Annetta Aramova, J\"urgen Herzog "p-Borel principal ideals", Illinois J.Math.41,no 1.(1997),103-121.
\bibitem[2]{ert}D.Eisenbud, A.Reeves, B.Totaro "Initial ideals, veronese subrings and rates of algebras", Adv.Math. 109 (1994), 168-187.
  \bibitem[3]{hpv}J\"urgen Herzog, Dorin Popescu, Marius Vladoiu "On the Ext-Modules of ideals of Borel
  type", Contemporary Math. 331 (2003), 171-186.
  \bibitem[4]{hp}J\"urgen Herzog, Dorin Popescu "On the regularity of p-Borel ideals", Proceed.of AMS, Volume 129,
             no.9, 2563-2570.
\bibitem[5]{par}Keith Pardue, "Non standard Borel fixed ideals", Dissertation, Brandeis University, 1994.
 \bibitem[6]{p}Dorin Popescu "Extremal Betti numbers and regularity of Borel type
 ideals", Bull. Math. Soc. Sc. Math. Roum. 48(96), no 1, (2005),
 65-72.
\end{thebibliography}
\end{document}